\documentclass[runningheads]{llncs}

\usepackage[T1]{fontenc}
\usepackage{graphicx}

\usepackage{amsfonts}
\usepackage{amsmath}
\usepackage{bm}
\usepackage{booktabs}
\usepackage{cancel}
\usepackage{mathtools}
\usepackage{siunitx}
\usepackage[table]{xcolor}
\usepackage[colorlinks=True,linkcolor=blue,citecolor=green,urlcolor=red]{hyperref}

\begin{document}
\title{Disease Progression Modelling and Stratification for detecting sub-trajectories in the natural history of pathologies: application to Parkinson's Disease trajectory modelling}
\titlerunning{DP-MoSt for detecting sub-trajectories in the natural history of pathologies}
\author{Alessandro Viani\inst{1} \and
Boris A. Gutman \inst{2} \and
Emile d'Angremont\inst{3} \and Marco Lorenzi\inst{1}}
%
\authorrunning{A. Viani B. Gutman E. d'Angremont M. Lorenzi}
%
\institute{Epione Research Team, Inria Center of Universit\'e C\^ote d’Azur, Sophia Antipolis, France, \and
Illinois Institute of Technology, Department of Biomedical Engineering, USA\and
Amsterdam University Medical Center, Department of Anatomy and Neurosciences, The Netherlands
\\
\email{alessandro.viani@inria.fr}}
\maketitle          	

\begin{abstract}
Modelling the progression of Degenerative Diseases (DD) is essential for detection, prevention, and treatment, yet it remains challenging due to the heterogeneity in disease trajectories among individuals. Factors such as demographics, genetic conditions, and lifestyle contribute to diverse phenotypical manifestations, necessitating patient stratification based on these variations. Recent methods like Subtype and Stage Inference (SuStaIn) have advanced unsupervised stratification of disease trajectories, but they face potential limitations in robustness, interpretability, and temporal granularity.
To address these challenges, we introduce Disease Progression Modelling and Stratification (DP-MoSt), a novel probabilistic method that optimises clusters of continuous trajectories over a long-term disease time-axis while estimating the confidence of trajectory sub-types for each biomarker.

We validate DP-MoSt using both synthetic and real-world data from the Parkinson’s Progression Markers Initiative (PPMI). Our results demonstrate that DP-MoSt effectively identifies both sub-trajectories and sub-populations, and is a promising alternative to current state-of-the-art models.

\end{abstract}

\section{Introduction}
\label{sec:introduction}
Modelling the progression of Degenerative Diseases (DD) is crucial for detection, prevention, and treatment purposes \cite{mould2012models}. This task is challenging due the generally large heterogeneity of disease trajectories  observed across subjects.
Despite a common degenerative process, the manifestation of symptoms and the configuration of biomarkers may vary widely among individuals, due for example to demographics, genetics conditions and life-style \cite{yang2011quantifying}. For this reason researchers are steering their attention towards the problem of \textit{patients stratification} based on  phenotypical manifestations of the disorder \cite{greenland2019clinical}.

Typical disease progression modelling approaches generally focus on estimating long-term biomarkers evolution from short-term patients observation. For instance, the Gaussian Process Progression Model (GPPM) \cite{lorenzi2019probabilistic} and the Personalized Input-Output Hidden Markov Model \cite{severson2020personalized} characterise the transition of biomarkers over time from normal to pathological stages, based on the assumption of an underlying disease trajectory defined by an absolute time axis. On the other hand, other approaches focus on the detection of sub-populations, for example biomarkers values \cite{zhang2016bayesian} or genetic observations \cite{whitwell2012neuroimaging}.

In the last decade, innovative methodologies have attempted to automatically stratify sub-types of disease trajectories. 
SuStaIn is a popular unsupervised method detecting sub-populations and their respective trajectories within a given dataset of patients and control population \cite{young2018uncovering,young2021ordinal,young2023subtype}. SuStaIn has been demonstrated in a variety of applications, showing its ability to identify sub-populations exhibiting common patterns of biomarkers changes \cite{aksman2021pysustain,zhou2023two}. Nevertheless, from an analytical perspective, SuStaIn presents some relevant limitations in terms of robustness and interpretability. First, disease progression is described as a discrete series of events, which limits the interpretability and temporal granularity of the estimated trajectories. Second, SuStaIn assumes the existence of cutoff values optimising separation between disease stages across biomarkers; the assumptions behind this statistical construct (e.g.\;Gaussian) are non-necessarily realistic about the biomarkers distribution across stages, and may negatively affect the robustness of the clustering task. Third, SuStaIn does not directly quantify the uncertainty in whether each biomarker exhibits a distinct trajectory between subtypes. For example, although certain biomarkers may not be discriminating between sub-types, SuStaIn will attempt at estimating group-specific cutoff values which may lead to poor interpretability and generalisation of the results. Finally, SuStain is designed to model cross-sectional information, without accounting for the temporally correlated nature of patients' time-series. Recent extensions of SuStain (t-SuStaIn) \cite{young2023subtype}, attempt to overcome this latter limitation by accounting for multiple measures per subjects. However, this approach still presents the above mentioned limitations, as it relies on the estimation of sequences of events occurring in a discrete space.

Disease Course Mapping (DCM) \cite{schiratti2017bayesian} is an orthogonal approach to SuStain, in which continuous parametric disease trajectories are optimized by accounting for patient's random effects represented by time-warp functions. DCM was recently extended to account for mixture of trajectories (MM-DCM) \cite{poulet2021mixture}. However, the mixture model there proposed assumes that all biomarkers' trajectories should be split into sub-progressions. This approach thus does not allow uncertainty quantification of  the split across biomarkers, and ultimately does not account for the specificity of certain biomarkers in characterising disease sub-types. 

To address these limitations, in this paper we present Disease Progression Modelling and Stratification (DP-MoSt), a novel probabilistic method to identify differential disease progression trajectories in heterogeneous cohorts.  DP-MoSt relies on the optimisation of clusters of continuous trajectories across a long-term disease time-axis, while also estimating the confidence for the existence of trajectories sub-types for each biomarker.

We validate our model on both synthetic and real-world data from the Parkinson’s Progression Markers Initiative (PPMI) \cite{aleksovski2018disease}. Our results show that our model is a promising alternative to the state-of-the-art, effectively identifying sub-trajectories and sub-populations, while providing interpretable and clinically meaningful solutions.

The manuscript is structured as follows: in Section \ref{sec:model_definition} we introduce DP-MoSt, and in Section \ref{sec:eem_optimisation} we describe the optimisation procedure and the statistical assumptions. In Section \ref{sec:toy_example} we present results on a panel of synthetic benchmarks, and in Section \ref{sec:ppmi_example} we provide a comparison between DP-MoSt and SuStaIn on clinical data from the PPMI database.

\section{Disease Progression Modelling and Stratification model}
\label{sec:proposed_method}
In this section, we provide the mathematical details of DP-MoSt, by describing the underlying statistical assumptions along with the associated optimization procedure.



\subsection{Model definition}
\label{sec:model_definition}

DP-MoSt is based on the optimization of two complementary problems: (i) estimating an absolute long-term disease time axis from short-term observations, and (ii) identifying along this disease time axis the existence of sub-populations with respective sub-trajectories.

Considering problem (i), for each individual $j$ we define the observations across all biomarkers as $\bm x^j = (\bm x^j_b)_{b=1}^B$; where $\bm x^j_b = (\bm x_b^j(\tilde t_{1}), \ldots, x_{1}^{j}(\tilde t_{k_j}))$ and $B$ is the number of biomarkers. Without loss of generality, for notational convenience, we assume that the sampling times are common among all subjects and biomarkers. To map the individual observations to a common disease time scale, we parameterize the individual time axis via a translation by a time-shift $\delta {\tilde t}_j$, i.e.\;$\bm t_j = \tilde t_{1:k_j} + \delta {\tilde t}_j$.
In this work, we evaluate the time shifts relying on the Gaussian process theory of GPPM \cite{lorenzi2019probabilistic}, which is based on the monotonic description of biomarkers trajectories from normal to pathological stages.

Considering problem (ii), given the measured observations $\bm x= \bm x^{1:J}$, where $J$ is the number of subjects, and the estimated absolute time $\bm t= \bm t_{1:J}$, we define a trajectory mixture model to identify the existence of sub-populations.

To achieve our goal, we assume that the evolution of each biomarker $b$ can be split into multiple sub-trajectories with probability $\xi_b$ (${b=1}, \ldots, B$). Once a sub-trajectory is considered, we assume that each subject $j$ is issued from this trajectory with probability $\pi_j$ (${j=1}, \ldots, J$). We observe that both $\xi = (\xi_b)_{1}^B$ and $\pi = (\pi_j)_1^J$ are independent with respect to time. This allows our model to link information deriving from longitudinal data; we also note that the probability for a subject belonging to one sub-trajectory must be consistent across all the biomarkers.

To simplify the inference process, compatibly with the monotonic assumption of GPPM, we adopt a parametric approach for the disease trajectories assuming that biomarkers follow increasing sigmoidal functions 
over time. 
We furthermore assume that the given measures are perturbed by additive Gaussian noise with standard deviation $\sigma=(\sigma_{b})_{b=1}^B$. Given the assumptions above, the posterior distribution for our model can be written as:
\begin{equation}
\begin{split}
	p(\theta, \sigma, \xi, \pi \mid \bm x) \propto p(\theta, \sigma, \xi, \pi)\prod_{j,b} p(\bm x_b^j \mid \theta_b, \sigma_b, \xi_b, \pi_j)
	\end{split}
	\label{eq:post_full}
\end{equation}
where for simplicity we omitted the conditioning on the time points. We observe that Equation \eqref{eq:post_full} implicitly assumes independence between the unknown parameters as well as independence between different subjects and biomarkers.

We can rewrite the equation by expanding the likelihood function in order to highlight the two-level mixture model formulation. In this setting, a first level deals with the sub-trajectory discovery task, while a second one determines the probability of a subject to belong to the sub-trajectory:
\begin{equation}
\begin{split}
	p(\bm x \mid \theta, \sigma, \xi, \pi)=&\prod_{j,b}\biggl[p(\bm x_b^j \mid \theta_b^0, \sigma_b)\xi_b +\\& \left(\pi_jp(\bm x_b^j \mid \theta_b^1, \sigma_b) + (1-\pi_j) p(\bm x_b^j \mid \theta_b^2, \sigma_b) \right)(1-\xi_b)\biggl],
	\end{split}
	\label{eq:like_full}
\end{equation}
where $p(\bm x_b^j \mid \theta_b^i, \sigma_b) = \prod_{\ell=1}^{k_j} \text{NormPDF}\left(x_b^j(t_\ell), f(t_\ell\mid \theta_b^i), \sigma_b\right)$ due to the assumption of additive Gaussian noise, and $f(t_\ell\mid \theta_b^i)$ is a Sigmoid function with parameters $\theta_b^i$ evaluated at $t_\ell$.


\subsection{Two-levels Expectation-Maximization}
\label{sec:eem_optimisation}
Model \eqref{eq:post_full} accounts for a substantial number of parameters. For each biomarker it includes: three parameters for each of the three Sigmoid functions $\theta_b^{0:2}$; the noise parameters $\sigma_b$; the probability for the existence of sub-trajectories $\xi_b$; and the membership probabilities for each subject $\pi$. To limit the computational cost, we focus on Maximum a Posteriori (MAP) estimation through Expectation Maximisation (EM), exploiting the two-levels mixture nature of the model:
\begin{equation}
\begin{split}
\hat\theta, \hat\xi, \hat\pi, \hat\sigma&= \arg\max\;\ln(p(\theta,\sigma,\xi,\pi\mid\bm x)) \\&= \arg\max\;\ln(p(\bm x\mid \theta,\sigma,\xi,\pi)) + \beta \sum_b \xi_b - \beta_N \sum_b \left(\ln(\sigma_b)+\frac{1}{\sigma_b}\right).
   \label{eq:loss}
\end{split}
\end{equation}

We rely on the following prior assumptions:
\begin{itemize}
	\item $p(\theta)\propto1$ and $p(\pi)\propto1$ as improper uniform prior for the Sigmoid parameters and sub-population subdivision in order to encode lack of information;
	\item $p(\sigma)=\prod_bp(\sigma_b)$, where $p(\sigma_b)$ is an Inv-Gamma distribution of shape parameter $\beta_N-1$ and scale parameter $\beta_N>1$ as classical prior distribution for the noise standard deviation \cite{calvetti2009conditionally} in order to penalise small values;
	\item $p(\xi)=\prod_bp(\xi_b)$, where $p(\xi_b)$ is a Laplace distribution \cite{kaban2007bayesian} of location $1$ and scale parameter $1/\beta$ restricted to the interval $[0,1]$, as regularisation term for penalising the introduction of a sub-trajectory to prevent overfitting.
\end{itemize}

DP-MoSt is designed to provide interpretable progression dynamics:  first, we estimate the parameters for the continuous biomarkers trajectories; second, we estimate the probability for each biomarker to present sub-trajectories; finally, we estimate for each subject its probability to belong to each sub-population. 

\section{Results}
\label{sec:results}
In this section, we validate the model on two different experimental scenarios\footnote{The code for DP-MoSt can be found at \url{https://github.com/alessandro-viani/DP-MoSt.git}.}: first, we demonstrate the effectiveness of the method in identifying sub-trajectories and sub-populations on an extensive synthetic benchmark; second, we apply the model to the real data from the PPMI dataset, comparing the solution to state-of-the-art solutions as provided by SuStaIn. The decision to validate only the clustering component of the model stems from the fact that the time-shift is assessed using GPPM, whose effectiveness has already been proven.

\subsection{Experiment on synthetic data}
\label{sec:toy_example}
\subsubsection{Data generation.}
We analyze the performance of DP-MoSt on the task of sub-trajectories and sub-populations identification by evaluating its performances across increasing levels of data complexity. This involves altering both the number of biomarkers and the Signal to Noise Ratio (SNR) \cite{johnson2006signal} between sub-trajectories. Specifically, we consider three different sets of biomarkers with $ B = 2, 5, 10 $ and assess performance under three distinct SNR conditions: low, normal, and high. This systematic variation allows us to test the robustness and accuracy of our model across a range of realistic scenarios of increasing difficulty.

For each configuration, we generate 100 datasets with the following parameters:

\begin{itemize}
	\item $J=100$ subjects
	\item $k_{1:J}=1$, i.e.\;one time point for each subject, randomly sampled in $[0, 20]$.
	\item $\sigma_b = 0.5$ to ensure a reasonable amount of noise in the data.
	\item Three fixed thresholds (low = 0.1, normal = 0.5, and high = 1) for the Mean Squared Error (MSE) between sub-trajectories to achieve a controlled SNR.
	\item Half of the biomarkers exhibit sub-trajectories (in the case of 5 biomarkers, 3 of them show the split).
	\item Equal partitioning of subjects between sub-populations.
	\item Parameters for the Sigmoid function describing the biomarker trajectories are randomly sampled from Gaussian distributions with parameters ensuring a positive supremum and rate of growth.
\end{itemize}

\subsubsection{Model setup.}
For each simulated dataset, we evaluate the parameters of the model described by Equation \eqref{eq:post_full} using the EM method described in Section \ref{sec:eem_optimisation}.

We initialize the model as follows to avoid local maxima for the posterior distribution:
\begin{itemize}
	\item set \(\xi_b = 0.5\) and \(\pi_j = 0.5\), assuming a complete lack of information on the sub-trajectories and sub-population probabilities;
	\item Initialize the noise standard deviation as the standard deviation of the data;
	\item randomly select the parameters for the Sigmoid functions from Gaussian distributions ensuring to provide a positive rate of growth and supremum;
	\item Set the prior parameters \(\beta\) and \(\beta_N\) as the 15\% of the number of subjects; this choice ensures effective regularisation while maintaining the values within a reasonable range.
\end{itemize}

\subsubsection{Performance metrics.}
To validate the proposed method, we use different metrics to evaluate the error in parameter approximations. For evaluating the error on the biomarker trajectories, we employ the Optimal Subpattern Assignment (OSPA) metric \cite{4567674,viani2020bayes}; this metric is particularly suitable because it accounts for potential differences in the number of true and estimated configurations. The OSPA metric measures the minimum MSE between the actual configuration of Sigmoids and the estimated one:
\begin{equation}
	\text{OSPA}(\theta, \hat \theta) = \min_{\phi}\sum_{i=1}^{\min\{\hat{d}, d\}} \frac{1}{N}\sum_j^N\left\|f(t_j\mid\hat\theta_{i}) - f(t_j\mid\theta_{\phi(i)})\right\|_2^2
	\label{eq:ospa}
\end{equation}
where $\hat \theta$ and $\hat d$ represent the estimated parameters and number of sub-trajectories, respectively; $\theta$ and $d$ denote the true values of the parameters; and the symbol $\phi$ represents all possible permutations.


\begin{figure}[t!]
	\centering
	\includegraphics[width=\textwidth]{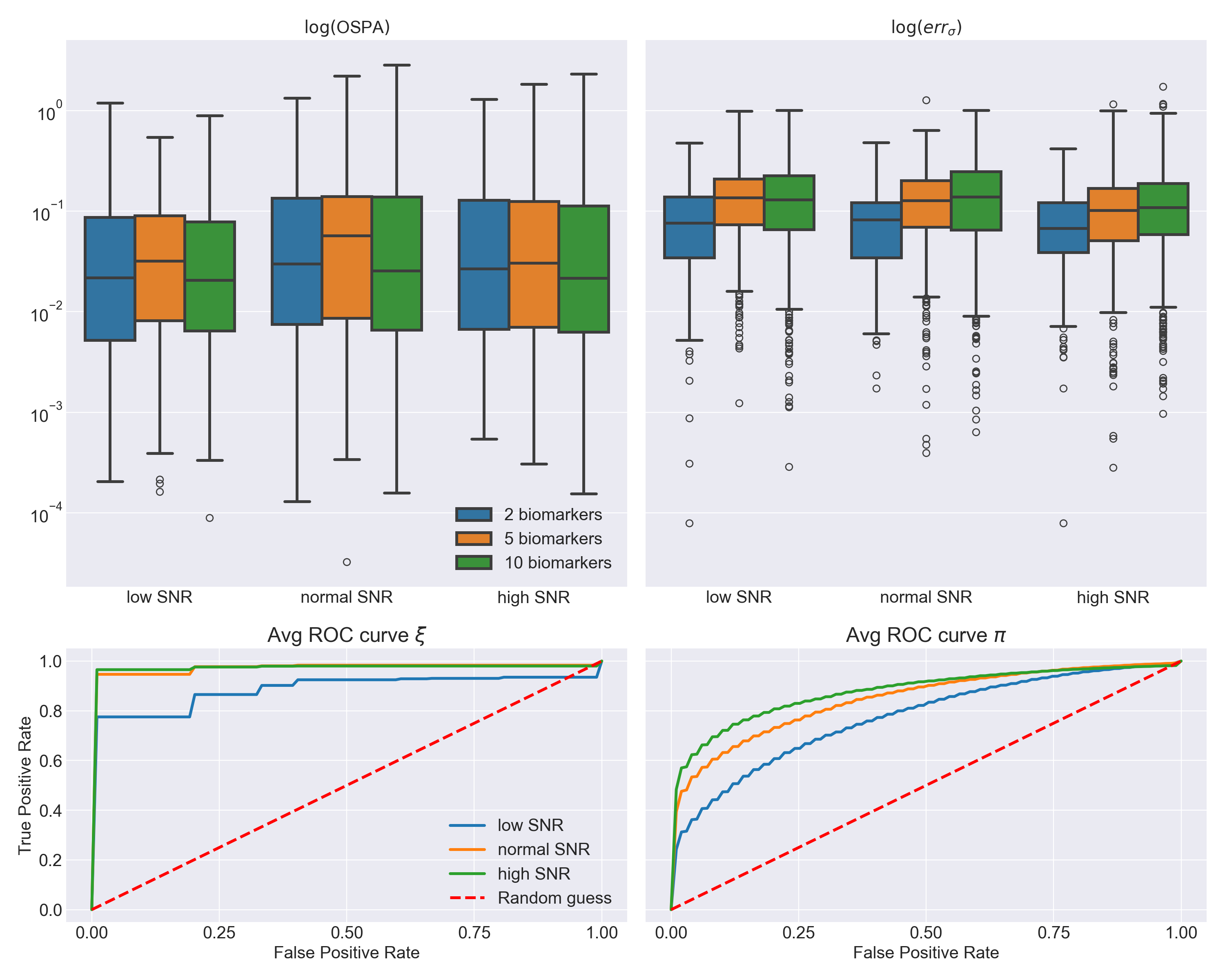}
	\caption{The figure shows performance metrics obtained on synthetic data. In the first row we show the OSPA error for the trajectories approximation and the error on the noise standard deviation; in the second row we show the ROC curves respectively for trajectory and individual clustering, as identified by parameters $\xi$ and $\pi$.}
	\label{fig:errors}
\end{figure}
\subsubsection{Results.}
The obtained results are summarised in Figure \ref{fig:errors}, where we show:
\begin{itemize}
\item First Panel, First Row: the logarithm of the OSPA error subdivided by the number of features considered and SNR levels. We observe that variations in data complexity do not significantly impact the trajectory approximation; indeed the error remains roughly constant over different data configurations.
\item Second Panel, First Row: the logarithm of the relative error on the noise standard deviation subdivided by the number of features considered and SNR levels. We can observe that the approximation error for the noise std increases with the data complexity.
\item First/Second Panel, Second Row: the ROC curve for the parameters $\xi$ and $\pi$
. We observe that the method performs better in estimating the number of sub-trajectories compared to the estimation of sub-populations, probably due to the fewer number of parameters to be estimated.
\end{itemize}

\subsection{PPMI data}
\label{sec:ppmi_example}
In this section we provide the results obtained considering the PPMI dataset, a comprehensive, multi-center longitudinal study for Parkinson's research.

The PPMI dataset includes extensive clinical, imaging, and biological data collected from PD patients, individuals with PD risk factors, and healthy controls. 
The overall data considered in this work is composed of a total of 3559 patients and 93 measured biomarkers.
Each patient is labelled according to the following clinical scores: tremor-dominant (TD), postural instability gait disorder predominant (PIGD), and intermediate \cite{aleksovski2018disease,potter2022serum}. 

We focus our analysis on a specific set of measures:
\begin{itemize}
	\item \textit{PIGD\_score}: an indicator that asserts the gravity of the PIGD classification;
	\item \textit{\textit{TD\_score}}: an indicator asserts the gravity of the TD classification;
	\item \textit{MCATOT}: the total amount of Montreal Cognitive Assessment scores \cite{julayanont2017montreal};
	\item \textit{NP2PTOT} \& \textit{NP3PTOT}: are the total score of MDS-UPDRS part 2 and 3. 
\end{itemize}
After preprocessing the data composed by the 5 variables of interest, we obtain a dataset including  \num{10198} longitudinal data from \num{1954} patients, with an average of 5 time points for each subject.
Data is analysed by applying DP-MoSt and SuStain. When applying DP-MoSt, we did not account for the temporal correlation between observations across subjects. This allows for a more fair comparison with SuStAin, where data are treated disregarding temporal dependency.

\begin{figure}[t!]
	\centering
	\includegraphics[width=\textwidth]{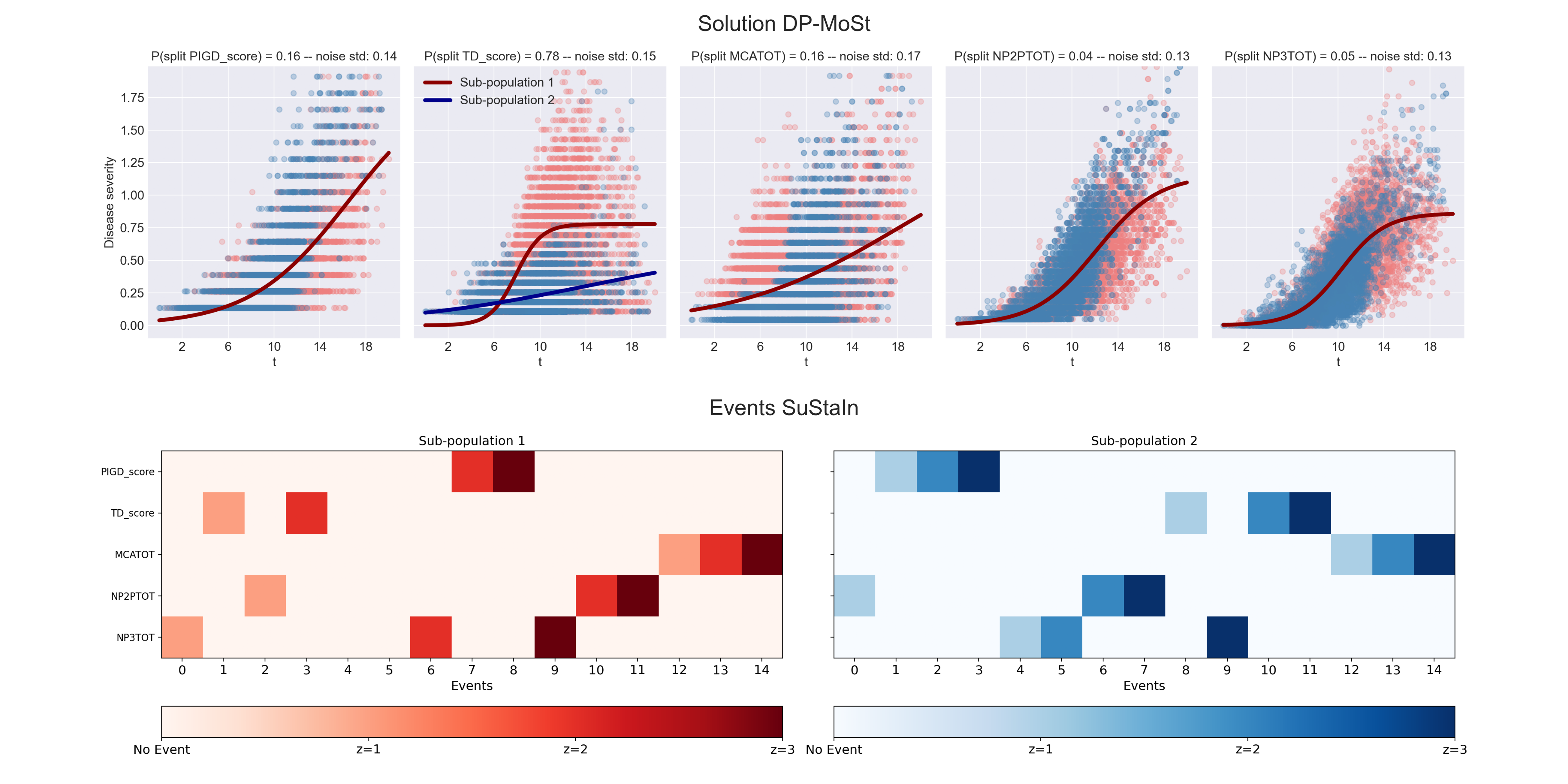}
	\caption{The figure shows the results obtained on the PPMI dataset. The first row shows the results obtained with DP-MoSt: (sub-)trajectories are represented as solid coloured lines, and subjects are colour-coded based on their estimated subgroups. The second row shows the event progression estimated by SuStaIn as two coloured matrices, one for each sub-population.}
	\label{fig:ppmi_results}
\end{figure}

\subsubsection{SuStaIn setup.}
SuStaIn estimates the maximum likelihood (ML) solution for the number of sub-types in the dataset as well as the associated sequence of events (i.e., stages of severity increase). Based on previous works, we decided to implement the z-score SuStaIn with three different progression stages associated to the quartiles $(z_1,z_2,z_3)$ and \num{100000} MCMC steps.

\subsubsection{DP-MoSt setup.}
We assess DP-MoSt by considering the same parameter values as for the experiment on synthetic data (Section \ref{sec:toy_example}). However, differently from the synthetic experiment in which we evaluated solely the clustering step of our method, we also optimise individual time-shift parameters, establishing an absolute time axis for the disease progression.

\begin{figure}
	\centering
    \includegraphics[width=\textwidth]{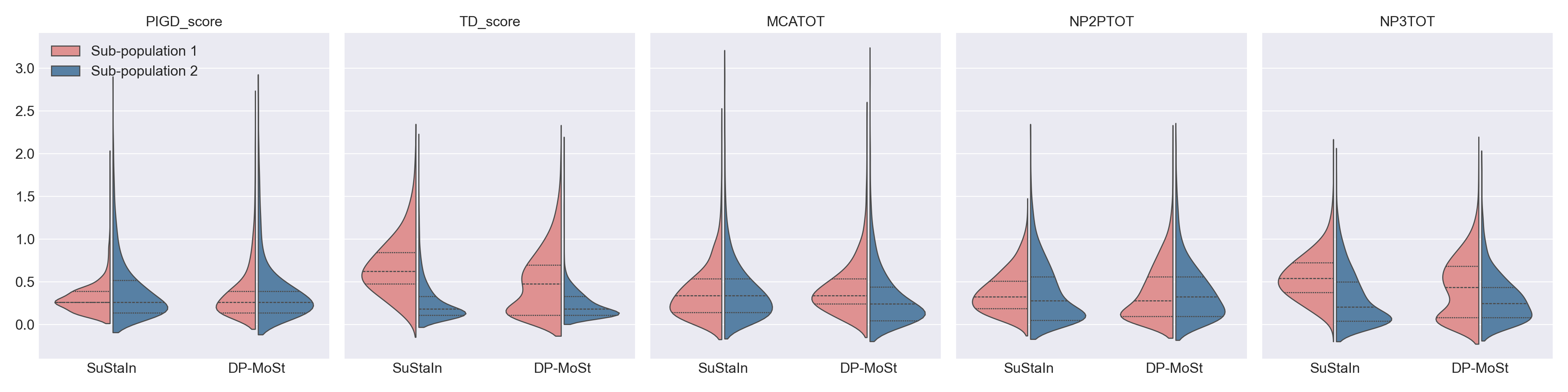}
	\caption{The figure shows the distribution for the biomarker values across the two different sub-populations detected by SuStaIn and DP-MoSt, where  we highlight the quartiles of the distributions as dotted lines.}
	\label{fig:ppmi_distribution}
\end{figure}

\subsubsection{Results.}
In Figure \ref{fig:ppmi_results} (first row) we show the solution for DP-MoSt, representing biomarker trajectories as solid lines with respect to the absolute time axis for disease progression (in years), and subjects as dots colour-coded based on their estimated sub-group. We observe that there is one clinical score that primarily governs the split into two sub-groups ($\xi_b>0.5$): the \textit{TD\_score}. The two sub-groups that result from this split show a clear subdivision with respect to \textit{TD\_score} values: sub-population 1 is associated with higher values for this biomarker compared to sub-population 2, with a clear differentiation of trajectories after t=6 years on the reparameterized time axis. We  also observe that, in spite of the split identified for the \textit{TD\_score} trajectory, the progressions of the other clinical scores do not exhibit clear partitioning in sub-trajectories ($\xi_b<0.5$). This aspect is illustrated in Figure \ref{fig:ppmi_distribution} where, besides the \textit{TD\_score}, the clinical scores' distributions between sub-populations show a substantial overlap. We finally note that, due to the continuous nature of DP-MoSt trajectories,  the resulting clinical scores' distributions can be multi-modal, as they describe the entire disease history across pathological stages (e.g. \textit{TD\_score} in Sub-population 1). This is less evident for SuStAin (Figure \ref{fig:ppmi_distribution}, left panels), which enforces a stronger separation between distributions across sub-groups, probably due to the discrete assumption for the events ordering. 

In Table \ref{tab:subdivision} (first column), we show the distribution of subjects into different sub-populations based on their label available at the last visit (Intermediate, PIGD, TD). We can observe that the sub-trajectories identified by DP-MoSt are respectively associated to subjects labeled as TD and PIGD, providing a clinically meaningful partitioning of the subjects. More specifically, DP-MoSt associates to sub-population 1 the majority of TD subjects (76\%) and to sub-population 2 the majority of PIGD subjects (62\%). The subjects labelled as Intermediate are not clearly assigned into any specific sub-population, probably because of their mixed composition.


In Figure \ref{fig:ppmi_results} (second row) we show the solution identified by SuStaIn as two coloured matrices, one for each sub-population, indicating the event ordering for each biomarker.
SuStaIn identifies 2 sub-populations, with sub-population 1 characterised by a faster increase of the \textit{TD\_score} and a slower change of the \textit{PIGD\_score} compared to sub-population 2. Notably, we can observe that the \textit{MCATOT} score does not show any differences between sub-populations, suggesting that this biomarker is not relevant to disease sub-typing.

In Table \ref{tab:subdivision} (second column) we observe that the clustering obtained with SuStaIn identifies two sub-populations that are unbalanced, with the majority of subjects (72\%) in the first sub-population. Therefore, differently from DP-MoSt, the stratification provided by SuStaIn does not provide a clear subdivision between clinical groups, substantially including the majority of subjects into the first sub-population. 

\begin{center}
\begin{table}[t]
\centering
\begin{tabular}{ccccccc}
\toprule
\multicolumn{1}{c}{} & \multicolumn{2}{c}{\textbf{DP-MoSt}}& \multicolumn{2}{c}{\textbf{SuStaIn}} \\
\cmidrule(rl){2-3} \cmidrule(rl){4-5}\cmidrule(rl){6-7}
\textbf{Condition} & {Sub-pop 1} & {Sub-pop 2} & {Sub-pop 1} & {Sub-pop 2} \\
\midrule
Intermediate &\cellcolor{red!28}0.56&\cellcolor{blue!22}0.44&\cellcolor{red!42}0.84&\cellcolor{blue!8}0.16 \\
PIGD & \cellcolor{red!19}0.38& \cellcolor{blue!31}\textbf{0.62} &\cellcolor{red!38}0.78& \cellcolor{blue!11}0.22 \\
TD & \cellcolor{red!38}\textbf{0.76}&\cellcolor{blue!12}0.24&\cellcolor{red!24}0.48&\cellcolor{blue!26}0.52 \\
\midrule
N°\% data & \cellcolor{gray!56}56\% & \cellcolor{gray!44}44\% & \cellcolor{gray!72}72\% & \cellcolor{gray!28}28\% \\
\bottomrule
\end{tabular}
\vspace{0.1cm}
\caption{The Table shows the subdivision between different sub-populations considering  solution provided by DP-MoSt and SuStaIn.}
\label{tab:subdivision}
\end{table}
\end{center}
\section{Conclusions}
In this paper, we introduced DP-MoSt, an innovative probabilistic method for identifying continuous biomarker trajectories and stratifying sub-populations in the context of DPMs. DP-MoSt addresses key limitations of current state-of-the-art methods by ensuring interpretability and robustness through a two-level mixture model that captures sub-population clusters while incorporating temporal information from longitudinal data. Additionally, DP-MoSt provides uncertainty quantification for both sub-trajectories and sub-population composition, enabling the characterization of pathological trajectory patterns.

We validated the model's performance with synthetic data, demonstrating its effectiveness. Furthermore, its application to the PPMI dataset yielded  interpretable and clinically relevant results; this feature can be due to the continuous formulation of DP-MoSt, improving the model's reliability and interpretability.

It is important to observe that SuStaIn provides an estimation of the number of sub-populations, a feature currently not included in DP-MoSt. However, we note that the introduction of multiple sub-populations with their respective sub-trajectories may also be obtained with DP-MoSt: this primarily involves more complex notation and higher computational costs, and is one of the improvements that will be implemented in future work.

Overall, DP-MoSt's ability to capture detailed sub-population characteristics makes it a promising tool for analysing heterogeneous disease progression patterns in longitudinal studies.

\begin{credits}
\subsubsection{\ackname}  This work has
been supported by the Michael J. Fox Foundation for Parkinson's Research (MJFF), and to the French government, through the 3IA Côte
d’Azur Investments in the Future project managed by the National
Research Agency (ANR) with the reference number ANR-19-P3IA-
0002, by the TRAIN project ANR-22-FAI1-0003-02, and by the ANR
JCJC project Fed-BioMed 19-CE45-0006-01. Data used in the preparation of this article were obtained on 05/01/2024 from the Parkinson’s Progression Markers Initiative (PPMI) database.

\subsubsection{\discintname}
The authors have no competing interests to declare that are
relevant to the content of this article.
\end{credits}
%
%

\bibliographystyle{splncs04}
\bibliography{bib}

\end{document}